\documentclass[leqno,12pt]{article} 
\usepackage{epic} 
\usepackage{ amsmath,amstext,amsthm,amssymb} 
\numberwithin{equation}{section} 

\begingroup   %  forthe spacing of the paragraphs 
\allowdisplaybreaks 
\swapnumbers 
\theoremstyle{plain}

\def\bthm#1.#2 #3\ethm{ 
\begin{\ifx#1ttheorem\fi% 
\ifx#1llemma\fi% 
\ifx#1ccorollary\fi% 
\ifx#1pproposition\fi% 
\ifx#1ddefinition\fi} \label{#1.#2} 
#3 \end{\ifx#1ttheorem\fi% 
\ifx#1llemma\fi% 
\ifx#1ccorollary\fi% 
\ifx#1pproposition\fi% 
\ifx#1ddefinition\fi}} 
 
\def\t#1/{theorem~\ref{t#1}}   \def\T#1/{Theorem~\ref{t#1}} 
\def\c#1/{corollary~\ref{c#1}}   \def\C#1/{Corollary~\ref{c#1}} 
\def\l#1/{lemma~\ref{l#1}}        \def\L#1/{Lemma~\ref{l#1}} 
\def\s#1/{section~\ref{s#1}} 
\def\e#1/{(\ref{e#1})} 
\def\d#1/{definition~\ref{d#1}} 
\def\f#1/{figure~\ref{f#1}} 
\def\Label #1 {\label{#1}}

\def\norm#1.#2.{\lVert#1\rVert_{#2}} 
\def\Norm#1.#2.{\bigl\lVert#1\bigr\rVert_{#2}} 
\def\NOrm#1.#2.{\Bigl\lVert#1\Bigr\rVert_{#2}} 
\def\NORm#1.#2.{\biggl\lVert#1\biggr\rVert_{#2}} 
\def\NORM#1.#2.{\Biggl\lVert#1\Biggr\rVert_{#2}} 
 
\def\ip#1,#2.{\langle #1,#2\rangle} 
\def\Ip#1,#2.{\bigl\langle#1,#2\bigr\rangle} 
\def\IP#1,#2.{\Bigl\langle#1,#2\Bigr\rangle}

\def\abs#1{\lvert#1\rvert} 
\def\Abs#1{\bigl\lvert#1\bigr\rvert} 
\def\ABs#1{\Bigl\lvert#1\Bigr\rvert} 
\def\ABS#1{\biggl\lvert#1\biggr\rvert}

\newcommand{\zb}{\ensuremath{\beta}}

\newcommand{\zd}{\ensuremath{\delta}} 
 
\newcommand{\ze}{\ensuremath{\epsilon}}

\newcommand{\zg}{\ensuremath{\gamma}}

\newcommand{\zI}{\ensuremath{\infty}} 
\newcommand{\zl}{\ensuremath{\lambda}}

\newcommand{\zm}{\ensuremath{\mu}} 
\newcommand{\zn}{\ensuremath{\nu}}

\newcommand{\zx}{\ensuremath{\xi}} 
 
\newcommand{\zs}{\ensuremath{\sigma}}

\newcommand{\zp}{\ensuremath{\pi}}
\newcommand{\zq}{\ensuremath{\chi}}

\def\z#1#2{\ifcase#1 {{\mathcal {#2}}}  %0=cal
\or {{\mathbf{#2}}}                    % 1=mathbold
\or { {\boldsymbol{#2}}}                  % 2=boldsymbol
\or{{\widetilde{#2}}}                   % 3=wide  tilde
\or {{\acute{#2}}} \or {{\grave{#2}}} \or {{\bar{#2}}} \or
{\dot{#2}} \or {\overline{#2}} \or {\underline{#2}}\fi}

\def\ZR{\ensuremath{\mathbb R}}
\def\ZD{\ensuremath{\mathbb D}}
\def\ZZ{\ensuremath{\mathbb Z}}

\def\ZT{\ensuremath{\mathbb T}}

\def\dist{\text{\rm dist}}
\def\bmo{{BMO}}
\def\seq{\preceq }   % for the funny less than
\def\mid{\,:\,}

\def\md#1#2\emd{\ifx0#1
\begin{equation*} #2 \end{equation*}\fi  %  single line display, no
\ifx1#1\begin{equation}#2\end{equation}\fi   % single line display,
\ifx2#1\begin{align*}#2\end{align*}\fi   % aligned display, no number
\ifx3#1\begin{align}#2\end{align}\fi    % aligned display, number
\ifx4#1\begin{gather*}#2\end{gather*}\fi  % multiline, not align, no
\ifx5#1\begin{gather}#2\end{gather}\fi   % multinline, not align
\ifx6#1\begin{multline*}#2\end{multline*}\fi  %  display too long for
\ifx7#1\begin{multline}#2\end{mutline}\fi  % as above, with numbers
}

\def\sqI{\sqrt{\abs I}}   

\def\ind#1{ {\mathbf 1}_{#1}}

\def\pp#1{P_{#1} }    %note the end of the argument is the first space
\def\ww #1 #2/{w_{#1}^{#2}}  %  For wavelets with super and
\def\w #1/{w_{#1}}      % as above, but with just subscript.

\def\commh{[M_b,H]}   

        \def\comm#1,#2.{[#1,#2]}

\def\mb{M_b}  % For the multiplication operator

\begin{document}

\title{A Characterization of Product $BMO$ by Commutators}
 \author{Michael T. Lacey\thanks{This work has been
  supported by an NSF grant, DMS--9706884.} \\
   Georgia Institute of Technology  \and Sarah H. Ferguson\thanks{This
work has been
  supported by an NSF grant, DMS--0071514}
\\   Wayne State University  }
\date{}

\maketitle

\section{Introduction}

 In this paper we establish a commutator estimate which allows one
to concretely identify the {\em product} $BMO$ space,
$ BMO(\ZR^{2}_{+} \times \ZR^{2}_{+})$, of Chang and R. Fefferman,
as an operator space on $L^2 (\ZR^2)$.
 The one parameter analogue of this result is a well--known theorem of Nehari
  \cite{n}.
The novelty of this paper is that we discuss a situation governed by a two
parameter family of dilations, and so the spaces $H^1$ and $BMO$ have a more 
complicated structure. 

 Here $\ZR^2_+$
denotes the upper half-plane and $BMO(\ZR^2_{+} \times \ZR^2_+)$
is defined to be the dual of the real-variable Hardy space
$H^1$ on the product domain $\ZR^{2}_+ \times \ZR^2_{+}$.
There are several equivalent ways to define this latter space and the reader
is referred to  \cite{gs} for the various characterizations.
We will be more interested in the
biholomorphic analogue of $H^1$, which can be defined in terms of the boundary
values of biholomorphic functions on  $\ZR^{2}_{+} \times \ZR^{2}_{+}$
and will be denoted throughout by $H^1( \ZR^{2}_{+} \times \ZR^{2}_{+})$, cf \cite{sw}.

In one variable, the space $L^2(\ZR)$ decomposes as
the direct sum $H^2(\ZR) \oplus \z8{H^2(\ZR)}$ where
$H^2(\ZR)$ is defined as the boundary values of functions in $H^2(\ZR^2_+)$
and $\z8{H^2(\ZR)}$ denotes the space of complex conjugate of functions in
$H^2(\ZR)$.  The space $L^2(\ZR^2)$, therefore, decomposes as the direct sum
of the four spaces $ H^2(\ZR) \otimes H^{2} (\ZR)$,
$\z8{H^2(\ZR)} \otimes H^{2} (\ZR)$, $ H^2(\ZR) \otimes \z8{H^{2}
(\ZR)}$ and $\z8{H^2 (\ZR)} \otimes \z8{H^2(\ZR)}$ where
the tensor products are the Hilbert space tensor products.
Let $P_{\pm,\pm}$ denote the orthogonal projection of
$L^2 (\ZR^2)$ onto the holomorphic/anti-holomorphic subspaces,
in the first and second variables, respectively, and let
$H_j$ denote the one--dimensional Hilbert transform in the
$j^{th}$ variable, $j=1,2$.   In terms of the projections $P_{\pm,\pm}$,
\[ H_{1} = P_{+,+} + P_{+,-} - P_{-,+} - P_{-,-} \quad \mbox{and} \quad
H_{2} = P_{+,+} + P_{-,+} - P_{+,-} - P_{-,-}.\]

The nested commutator determined
by the function $b$ is the operator $[[M_b,H_1],H_2]$ acting on
$L^2 (\ZR^2)$ where, for a function $b$ on the plane, we define $\mb f:=bf$.

In terms of the projections $P_{\pm,\pm}$,
it takes the form
\md1\Label e.commutator \tfrac 14 [[M_b,H_1],H_2]  =
{P_{+,+}} M_b P_{-,-}  - P_{+,-} M_b P_{-,+}  - P_{-,+} M_b
P_{+,-} + P_{-,-} M_b P_{+,+}. \emd

Ferguson and Sadosky \cite{fs} established the inequality
$    \norm [[M_b,H_1],H_2].L^2.\leq{}c\norm b.\bmo.$.  The
 main result is the converse inequality.

\bthm t.1 There is a constant $c > 0$ such that
$\norm b.\bmo. \leq{}c \norm [[M_b,H_1],H_2].L^2 \rightarrow
L^2.$ for all functions $b$ in $BMO(\ZR^2_{+} \times \ZR^2_+)$.
\ethm

As A.~Chang and R.~Fefferman have established for us, the structure of the 
space $BMO$ is  more complicated in the two parameter setting, 
requiring a more subtle approach to this theorem, despite the 
superficial similarity of the results to the one parameter setting. 
The proof  relies on three key ideas.  The first is
the dyadic characterization of the $BMO$ norm given in \cite{cf}.
The second is a variant of Journ\'e's lemma, \cite{journe}, (whose proof is included
in the appendix.)  The third idea is that we have the estimates,
the second of which was shown in \cite{fs},
\md0
\norm b .BMO(rect).{}\leq{}c
 \norm [[M_b, H_1],H_2].L^2 \rightarrow L^2.    {}\leq{}c'
  \norm b .BMO. .
\emd
An unpublished example of L.~Carleson shows that 
the rectangular $BMO$ norm is not comparable to the $BMO$
norm, \cite{f}.
We may assume that the rectangular $BMO$ norm of
the function $b$ is small. Indeed, this turns out to be an essential aspect 
of the argument.  

From \t.1/ we deduce a {\em{weak}} factorization
for the (biholomorphic) space $H^1 (\ZR^2_+ \times \ZR^2_+)$.
The idea is that if the function $b$ has biholomorphic extension to
$\ZR^2_+ \times \ZR^2_+$ then for functions $f,g \in L^2(\ZR^2)$,
\md0
 \tfrac 14 {\langle} [[M_b,H_1],H_2]f,g{\rangle}  =
{\langle} b, \overline{ P_{-,-}f}P_{+,+}g {\rangle}.
\emd
So in this case, the operator norm of the nested commutator
$[[M_b,H_1],H_2]$ is comparable to the {\em{dual}} norm
\md0
\norm b .\ast. := \sup \abs{{\langle}fg,b{\rangle}}
\emd
where the supremum above is over all pairs $f,g$ in the unit ball of
$H^{2} (\ZR^{2}_+ \times \ZR^{2}_+ )$.
On the other hand, since $\norm b.BMO.$ and
$\norm [[M_b,H_1],H_2].L^2 \rightarrow L^2.$ are comparable,
the dual norm above satisfies
\md0
\norm b . \ast. \sim \sup \abs{{\langle} h,b{\rangle}}
\emd
where the supremum is over all functions $h$ in the unit ball of
$H^1(\ZR^2_+ \times \ZR^2_+)$.

\bthm c.1 Let $h$ be in $H^1(\ZR^2_+ \times \ZR^2_+)$ with $\norm h .1. =1$.
Then there exists functions \( (f_j), (g_j) \subseteq H^2(\ZR^2_+ \times \ZR^2_+) \)
such that \( h = \sum_{j =1 }^{\infty} f_{j} g_{j} \) and
\( \sum_{j =1}^{\infty} \norm f_j .2. \norm g_j .2. \leq c. \) \ethm

We remark that the weak factorization above implies the analogous
factorization for $H^1$ of the bidisk.  Indeed, for all $1 \leq p <
\infty$, the map
\( u_{p} : H^p( \ZR^2_+ \times \ZR^2_+ ) \rightarrow H^p (\ZD^2)
\)
defined by
\[ (u_{p} f )(z,w) = \pi^{2/p} \left( \frac{2i}{1-z} \right)^{2/p}
\left( \frac{2i}{1-w} \right)^{2/p} f(\alpha (z), \alpha (w)) \quad
\alpha(\lambda) := i \frac{1 + \lambda}{1 - \lambda}, \]
is an isometry  with isometric inverse
\[ (u_{p}^{-1} g )(z,w) = \pi^{-2/p} \left( \frac{1}{z+i} \right)^{2/p}
\left( \frac{1}{w+i} \right)^{2/p} g(\beta(z), \beta (w)) \quad
\beta(\lambda) := \frac{\lambda - i}{\lambda + i}. \]

The dual formulation of weak factorization for $H^1 (\ZD^2)$ is a Nehari
theorem for the bidisk.  Specifically, if $b \in H^2(\ZD^2)$ then
the {\em little} Hankel operator with symbol $b$ is densely defined
on $H^2 (\ZD^2)$ by the formula
\md0
\Gamma_{b} f = P_{-,-} (\z8{b} f).
\emd
By \e.commutator/,
$\norm \Gamma_b . . = \norm [[M_{\z8{b}},H_1],H_2]. L^2 \rightarrow L^2.$
and thus, by \t.1/, $\norm \Gamma_b . .$ is comparable to
$\norm b. BMO. $ which, by definition, is just the norm of $b$ acting
on $H^1 (\ZD^2)$.  So the boundedness of the Hankel operator
$\Gamma_b$ implies that there is a function $\phi \in L^\infty (\ZT^2)$
such that $P_{+,+}\phi = b$.

Several variations and complements on these themes in the one parameter
setting have been obtained by Coifman, Rochberg and Weiss \cite{one}.

The paper is organized as follows.
Section 2 gives the one-dimensional preliminaries for the
proof of \t.1/ and  Section 3 is devoted to the proof of
\t.1/.  The appendix contains the variant
of Journ\'e's  lemma we use in our proof in Section 3.
One final remark about notation.
$A\seq B$ means that there is an absolute constant $C$ for which
$A\le{}CB$.  $A\approx B$ means that $A\seq B$ and $B\seq A$.

\medskip 

We are endebited to the anonymous referee.

\section{Remarks on the one dimensional case}

Several factors conspire to make the one dimensional case  easier
than the higher dimensional case. Before proceeding with the
higher dimensional case, we make several comments about the one
dimensional case, comments that extend and will be useful in the
subsequent section.

Let $H$ denote the Hilbert transform in one variable,
$\pp +  = \tfrac12(I+H)$ be
the projection of $L^2(\ZR)$ onto the positive frequencies, and
$\pp -$ is $\tfrac12(I-H)$, the projection onto the negative
frequencies.  We shall in particular rely upon the
following basic computation.
 \md1  \Label e.basic
\tfrac 12 \commh \overline b=\pp- \abs{\pp- b}^2-\pp+ \abs{\pp+
b}^2. \emd The frequency distribution of $\abs{\pp- b}^2$ is
symmetric since it is real valued.  Thus, \md2 \lVert
b\rVert_4^2\seq{} & \norm \pp- \abs{\pp- b}^2-\pp+ \abs{\pp+
b}^2.2.
\\
{}\le{}& \norm \commh.2\to2.\norm b.2. \emd Moreover, if $b$ is
supported on an interval $I$, we see that \md0 \norm b.2.\le{}
\abs{I}^{1/4}\norm b.4.\seq\abs{I}^{1/4} \lVert
\commh\rVert_{2\to2}^{1/2} \lVert b\rVert_2^{1/2} \emd which is
the $\bmo$ estimate on $I$.  We seek an extension of this
estimate in the two parameter setting.
\bigskip

We use a wavelet proof of \t.1/, and
specifically use a wavelet with compact frequency support
constructed by Y.~Meyer \cite{m}.  There is a Schwarz function $w$ with
these properties.
\begin{itemize}
\item $\norm w.2.=1$.
\item $\widehat w(\zx)$ is supported on   $[2/3,8/3]$ together with the
symmetric
interval about $0$.
\item $\pp\pm  w$ is a Schwartz function.  More particularly, we have
\md0 \abs{w(x)}, \abs{\pp \pm w(x)}\seq{}(1+\abs x)^{-n}, \qquad
n\ge1. \emd
\end{itemize}

Let $\z0D$ denote a collection of dyadic intervals on $\ZR$. For
any interval
  $I$, let $c(I)$ denote it's center, and define
  \md0
  \w I/(x):=\frac1{\sqI}w\Bigl(\frac{x-c(I)}{\abs I}\Bigr).
  \emd
  Set $\ww I \pm/:=\pp\pm \w I/$.  The central facts that we need
about
  the functions $\{\w I/\mid I\in\z0D\}$ are these.

  First, that these functions are an orthonormal basis on $L^2(\ZR)$.
Second,
  that we have the Littlewood--Paley inequalities, valid on all $L^p$,
though $p=4$
  will be of special significance for us.  These inequalities are
  \md1\Label e.LPaley
  \norm f.p.\approx \NORm \Bigl[\sum_{I\in\z0D}
  \frac{\abs{\langle f,\w I/\rangle}^2}{\abs I}\ind
I\Bigr]^{1/2}.p.,\qquad
  1<p<\zI.
  \emd
  Third,  that the functions $\w I/$ have good localization properties
in the
  spatial variables.  That is,
  \md1 \Label e.spatial
  \abs{\w I/(x)},\abs{\ww I \pm/(x)}\seq \abs{I}^{-1/2}
\zq_I(x)^n,\qquad n\ge1,
  \emd
  where $\zq_I(x):=(1+\dist(x,I)/\abs I)^{-1}$.  We find the compact 
localization of the wavelets in frequency to be very useful.  The price we pay 
for this utility below is the  careful accounting 
 of ``Schwartz tails" we shall make in the main argument. 
     Fifth, we have the
identity
  below for the commutator of one $\w I/$ with a $\w J/$.  Observe that
since
  $\pp + $ is one half of $I+H$, it suffices to replace $H$ by $\pp + $
in the
  definition of the commutator.
  \md3
  \nonumber
  \w I,J/:={}&\comm \w I/,\pp+ .\z8{\w J/}
  \\\nonumber{}={}&
  \w I/\z8{\ww J -/}-\pp +\w I/\z8{\w J/}
  \\\nonumber{}={}&
  \pp-\w I/\z8{\ww J -/}-\pp +\w I/\z8{\ww J +/}
  \\\nonumber{}={}&
  \pp-\ww I -/\z8{\ww J -/}-\pp +\ww I +/\z8{\ww J +/}
  \\{}={}&   \Label  e.wij
  \begin{cases}
  0 & \text{if $\abs I\ge4\abs J$.}   \\
  \pp- \abs{\ww I -/}^2-\pp + \abs{\ww I +/}^2& \text{if $I=J$.}  \\
  \ww I -/\z8{\ww J -/}-\ww I +/\z8{\ww J +/}& \text{if $\abs J\ge4\abs
I$.}
  \end{cases}
  \emd
From this we see a useful point concerning orthogonality.  For
intervals $I,I',J$ and $J'$, assume $ \abs J\ge8\abs I$ and
likewise for $I'$ and $J'$.  Then \md1 \Label e.orthoI
\text{supp}(\widehat{ \w I,J/})\cap\text{supp}(\widehat{\w
I',J'/})=\emptyset,\qquad \abs{I'}\ge8\abs I . \emd Indeed, this
follows from a direct calculation.  The positive frequency support
of $ \ww I +/\z8{\ww J +/}$ is contained in the interval $[(3\abs
I)^{-1},8  (3\abs I)^{-1}]$. Under the conditions on $I$ and $I'$,
the frequency supports are disjoint.

     %%%%%%%%%%%%%%%%%%%%  Proof of main Theorem
%%%%%%%%%%%%%%%%
 \section{Proof of the main theorem}

$BMO(\ZR^2_+ \times \ZR^2_+)$ will denote the $BMO$ of two
parameters (or {\em {product}} $BMO$) defined as the dual
of ({\em{real}}) $H^1(\ZR^2_+ \times \ZR^2_+)$.
The following characterization of the
space $BMO(\ZR^2_+ \times \ZR^2_+)$ is due to Chang and R. Fefferman \cite{cf}.

The relevant class of rectangles is $\z0R=\z0D\times \z0D$, all
rectangles which are products of dyadic intervals.
These  are indexed by $R\in\z0R$.  For such
a rectangle, write it as a product $R_1\times R_2$ and then define
\md0 v_{R}(x_1,x_2)=w_{R_1}(x_1)w_{R_2}(x_2). \emd
A function $f\in BMO(\ZR^2_+ \times \ZR^2_+)$ iff \md0
\sup_{U}\Biggl[\abs{U}^{-1}\sum_{R\subset U}\abs{\ip
f,v_R.}^2\Biggr]^{1/2}<\zI. \emd
Here, the sum extends over those
rectangles $R\in\z0R$ and the supremum is over all open sets in
the plane of finite measure.  Note that the supremum is taken over
a much broader class of sets than merely rectangles in the plane.
We denote this supremum as $\norm f.\bmo.$.  

In this definition,
if the supremum over $U$ is restricted to
just rectangles, this defines the ``rectangular $\bmo$" space, and
we denote this restricted supremum as $\norm f.\bmo(\text{rec}).$.

\medskip

Let us make a further comment on the $BMO$ condition.   Suppose
that for $R\in\z0R$, we have non--negative constants $a_R$ for
which \md0 \sum_{R\subset U}a_R\le{}\abs U, \emd for all open sets
$U$ in the plane of finite measure.  Then, we have the
John--Nirenberg inequality \md0 \NORm \sum_{R\subset U}
\abs{R}^{-1}a_R\ind R.p.\seq{}\abs U^{1/p},\qquad 1<p<\zI. \emd
See \cite{cf}.  This, with the Littlewood--Paley inequalities,
will be used several times below, and referred to as the
John--Nirenberg inequalities.

\subsection{The Principal Points in the Argument}

We begin the principle line of the argument. The function $b$ may be taken to be
of Schwarz class. By multiplying
$b$ by a constant, we can assume that the $BMO$ norm of $b$ is
$1$. Set $B_{2\to2}$ to be the operator norm of $[[M_b,H_1],H_2]$.
Our purpose is to provide a lower bound for $B_{2\to2}$.
Let $U$ be an open set of finite measure for which we have the
equality
\md0
\sum_{R\subset U}\abs{\ip b,v_R.}^2=\abs U.
\emd
As $b$ is of Schwarz class, such a set exists. 
By invariance under dilations by a factor of two, we can assume that
$\frac12\le\abs U\le1$.  In several estimates below, the measure
of $U$ enters in, a fact which we need not keep track of.

An essential point is that we may assume that the
rectangular $\bmo$ norm of $b$ is at most $\epsilon$.
The reason for this is that we have the estimate $ \norm
b.\bmo(\text{rec}).\seq{}B_{2\to2} $. See \cite{fs}.
Therefore, for a small constant $\epsilon$ to be chosen below, we can
assume that $\norm b.\bmo(\text{rec}).\seq\epsilon$, for otherwise we
have a lower bound on $B_{2\to2}$.

Associated to the set $U$ is the set $V$, defined below, which  is
an expansion of the set $U$.  
In defining this expansion, it is critical that the measure of $V$ be only
slightly larger than the measure of $U$, and so in particular we do not use the
strong maximal function to define this expansion. 
 In the definition of $V$, the parameter $\zd
> 0$ will be specified later and $M_j$ is the one dimensional
maximal function  applied in the direction $j$. Define 
 \md4
 V_{  ij}:=\{M_i\ind {\{M_j 
\ind{}U>\zd\}}>\zd\} ,  \\
  V:=V_{12}\cup{}V_{21} ,\\
\zm(R):=\sup\{\zm\mid \zm R\subset V\}\quad \text{$R\subset U$.}
\emd
 The quantity $\zm(R)$ measures how deeply a rectangle $R$ is
inside $V$. This quantity  enters into the essential Journ\'e's
Lemma, see \cite{journe} or the variant we prove in the appendix.

In the argument below, we will be projecting $b$ onto subspaces spanned by
collections of wavelets.  These wavelets are in turn indexed by
collections of rectangles.  Thus, for a collection
$\z0A \subseteq \z0R$, let us denote
\md4
b^{\z0A}:=\sum_{R\in\z0A}\ip b,v_{R}.v_{R}.
%\\
%\{\z0A,\z0B\}:={} [M_{b^{\z0A}},H_1],H_2]\z8{b^{\z0B}}={}
%\sum_{R\in\z0A}\sum_{{R'}\in\z0B} \ip b,v_{R}.\z8{\ip
%b,v_{R'}.}\{v_R,v_{R'}\}
%[M_{v_R },H_1],H_2]\z8{v_R}.
\emd
The relevant collections of rectangles are defined as
\md4 \z0U:=\{R\in\z0R\mid R\subset U\},
\\
\z0V=\{ R\in\z0R-\z0U\mid R\subset V\}, \\
\z0W=\z0R-\z0U-\z0V. \emd
 For functions $f$ and $g$, we
set $\{f,g \} := [ [ M_f , H_1 ], H_2 ] \z8{g}$.

We will demonstrate that  for all $\zd,\epsilon >0$ there is a
constant $K_{\zd} > 0$ so that
 \begin{itemize}
\item [(i)]$\norm
\{b^{\z0V},b^{\z0U}\}.2.\seq\zd^{1/8}$  
\item [(ii)] $\norm
\{b^{\z0W},b^{\z0U}\}.2.\leq{}K_{\zd}{}{\epsilon}^{1/3}$
\end{itemize}
Furthermore, we will show that $1\seq\norm \{b^{\z0U},b^{\z0U}\}.2.$.
Since
$b=b^{\z0U}+b^{\z0V}+b^{\z0W}$, $\norm b^{\z0U}.2.\seq 1$ and
$\zd,\epsilon>0$ are arbitrary, a lower bound on $B_{2\to2}$ will
follow from an appropriate choice of $\zd$ and $\epsilon$.
To be specific, one concludes the argument by estimating
\md2
1\seq{}&\norm \{b^{\z0U},b^{\z0U}\} .2.
\\{}\seq{}&\norm \{b^{\z0U}+b^{\z0V} ,b^{\z0U}\}.2.+ \zd^{1/8}
\\{}\seq{}&\norm \{b^{\z0U}+b^{\z0V}+b^\z0W ,b^{\z0U}\}.2.+
\zd^{1/8}+K_{\zd}{\epsilon}^{1/3}
\\{}\seq{}&B_{2\to2}+\zd^{1/8}+K_{\zd}{\epsilon}^{1/3}
\emd
Implied constants are absolute.  Choosing  $\zd$ first and
then $\epsilon$ appropriately small supplies a lower bound on
$B_{2\to2}$.

\bigskip         %%%%%%%%%%%%%%%%%  First Estimate

The estimate $1\seq\norm \{b^{\z0U},b^{\z0U}\}.2.$ relies on the
John--Nirenberg inequality and the two parameter version of
\e.basic/, namely 
\md0 \tfrac 14 [[M_b,H_1],H_2]\z8 b =
{P_{+,+}} \abs{P_{+,+} b}^2 - P_{+,-} \abs{P_{+,-}b}^2 - P_{-,+}
\abs{P_{-,+}b}^2 + P_{-,-} \abs{P_{-,-}b}^2 \emd
This identity easily follows from the one-variable identities. Here
$P_{\pm,\pm}$ denotes the projection onto the positive/negative
frequencies in the first and second variables. These projections
are orthogonal and moreover, since $\abs{P_{\pm,\pm} b}^2$ is real
valued we have that $\norm P_{\pm,\pm}\abs{P_{\pm,\pm}
b}^2.2.\ge\frac14\norm \abs{P_{\pm,\pm}b}^2.2.$.  Therefore,
$\lVert b^{\z0U}\rVert_4^2\seq\norm \{b^{\z0U},b^{\z0U}\}.2.$. It
follows by the John--Nirenberg inequality that \md2
1{}\seq{}&\norm b^{\z0U} .2.
\\= & \Biggl[ \sum_{R\in\z0U}\abs{\ip b,v_R.}^2\Biggr]^{1/2}
\\ \seq{}&
 \Biggl\lVert \Biggl[ \sum_{R\in\z0U}\abs{\ip b,v_R.}^2\ind
 R\Biggr]^{1/2}\Biggr\rVert_4^2
\\ \seq{}& \lVert  b^{\z0U} \rVert_4^2
\\ \seq{}&\norm \{b^{\z0U},b^{\z0U}\}.2.
 \emd

\bigskip  %%%%%%%%%%%%%%%%  Second estimate

 The estimate $(i)$ relies on  the fact that the one dimensional
maximal function maps $L^1$ into weak $L^1$ with norm one. Thus,
for all $0<\zd<1/2$,
\md0
\abs{ \{M_2\ind{ \{M_1\ind U>1-\zd\}  }>1-\zd\}}
\le(1-\zd )^{-2}\abs U\le{}(1+6\zd)\abs U.
 \emd
 Now, if
$R\in\z0V$, then $R\subset V$ and since $b$ has
$BMO$ norm one, it follows that
 \md0
  \abs U+\lVert
b^{\z0V}\rVert_2^2=\sum_{R\in\z0U\cup\z0V}\abs{\ip
b,v_R.}^2\le{}(1+6\zd)\abs U  .
\emd
 Hence $\norm
b^{\z0V}.2.\seq\zd^{1/2} $. Yet the $\bmo$ norm of $b^{\z0V}$
can be no more than that of $b$,
which is to say $1$.  Interpolating norms we see that $ \norm
b^{\z0V}.4.\seq\zd^{1/4}$, and so \md0 \norm
\{b^{\z0V},b^{\z0U}\}.2.\seq\norm b^{\z0V}.4.\norm
b^{\z0U}.4.\seq\zd^{1/4} . \emd

\bigskip

%%%%%%%%%%%%%%%%%%%%%%%%%%%   Third estimate
\subsection{Verifying the  Estimate $(ii)$}

We now turn to the estimate $(ii)$.
Roughly speaking $b^{\z0U}$ and $b^{\z0W}$ live on disjoint sets.
But in this argument we are trading off precise Fourier support
of the wavelets for imprecise spatial localization, that is the 
"Schwartz tails" problem. 
  Accounting
for this requires a careful
analysis, invoking several subcases.

A property of the commutator that we will rely upon is that it
controls the geometry of $R$ and $R'$.  Namely,
$\{v_{R'},v_{R}\}\not=0$ iff writing
$R=R_1\times R_2$ and likewise for $R'$, we have for both $j=1,2$,
$\abs{R_j'}\le{}4\abs{R_j} $. This follows immediately from our
one--dimensional calculations, in particular \e.wij/.  We
abbreviate this condition on $R$ and $R'$ as $R'\seq 
R$ and restrict our attention to this case.

Orthogonality also enters into the argument.  Observe the following.
For rectangles $ R^k,\z3R_k$, $k=1,2$, with $ \z3R^k\seq{}
  R^k$,  and for $j=1$ or $j=2$
\md1\label{e.doubleor}
\text{ if $8\abs{\z3R^1_j}\le{}\abs{R^1_j},\ \abs{\z3R^2_j}$
and
$8\abs{\z3R^2_j}<\abs{R^2_j}$ then  $\ip
 v_{\z3R^{1}}\z8{v_{R^1}},v_{\z3R^{2}}\z8{v_{R^2} }  .=0$.}
 \emd
 This follows from applying \e.orthoI/ in the $j$th coordinate.

Therefore, there are different partial orders on rectangles that are
relevant to our argument.  They are
\begin{itemize}
\item $R'<R$ iff $8\abs{R'_j}\le{}\abs{R_j}$ for $j=1$ and $j=2$.
\item  For $j=1$ or $j=2$, define $R'<_jR $ iff
$R'\seq R$ and $8 \abs{R'_j}\le{}\abs {R_j}$ but $R'\not<R$.
\item $R'\simeq R$ iff $\frac14\abs{R_j}\le\abs{R'_j}\le{}\abs{R_j}$ for $j=1$ and $j=2$.
\end{itemize}
These four partial orders divide the  collection $\{ (R',R)\mid R'\in\z0W,\ R\in\z0U,\
R'\seq R\}$ into four subclasses which require
different arguments.

In each of these four arguments, we have recourse to this
definition.   Set $\z0U_k$, for $k=0,1,2,\ldots$ to be
those rectangles in $\z0U$ with $2^{-k-1}<{}\zm(\z0U_k)\le2^k$.

Journ\'e's Lemma enters into the considerations.
Let $\z0U'\subset \z0U_k$ be a collection of
rectangles which are pairwise incomparable with respect to
inclusion.  For this
latter collection, we have the inequality
 \md1  \Label e.JU
 \sum_{ {R\in\z0U' }}\abs
R\leq{}K_{\zd}2^{k/100}\ABs{\bigcup_{R\in\z0U'}R}.
 \emd
 See Journ\'e \cite{journe}, also see the appendix. This
together with the assumption that $b$ has small rectangular $\bmo$ norm
gives us
 \md1\label{e.u8}
 \norm b^{\z0U_k}.\bmo.\le{}K_{\zd}2^{k/100}\epsilon.
 \emd
 This interplay between the small rectangular $BMO$ norm and
Journ\'e's Lemma is a decisive feature of the argument.

Essentially, the decomposition into the collections $\z0U_k$ is a
spatial decomposition of the collection $\z0U$.  A corresponding
decomposition of $\z0W$ enters in.  Yet the definition of this
class differs slightly depending on the partial order we are
considering.

\bigskip

 For $R'\in\z0W$ and $R\in\z0U$
the term $\{v^{R'},v^{R}\}$ is a linear combination of \md0
v_{R'}H_2H_1\z8{v_R},\quad H_2(v_{R'}H_1\z8{v_R}),\quad
(H_1v_{R'})(H_2\z8{v_R}),\quad H_1H_2 (v_{R'}\z8{v_R}).
\emd
Consider the last term.  As we are to estimate an $L^2$ norm, the
leading operators $H_1H_2$ can be ignored.  Moreover, the
essential properties of wavelets used below still hold for the
conjugates and Hilbert transforms of the same.  These properties
are   Fourier localization and    spatial localization.
Similar comments apply to the other three terms, and so the
arguments  below applies to each type of term seperately.

\subsubsection{The partial order `$<$'}
 
 We consider  the  case of $ R'  <R$ for $R'\in\z0W$ and 
$R\in\z0U$.
The sums we considering  are related to the following definition.
Set
\md0
b^{{\z0U_k}}_{\text{trun}}(x):=\sup_{R'}\ABS{\sum_{
\substack{R\in{\z0U_k}\\ R'<R}} \ip b,v_R. v_R(x)}.
\emd
Note that we consider the maximal truncation of the sum over all
choices of dimensions of the rectangles in the sum.  Thus, this
sum is closely related to the strong maximal function $M$ applied
to $b^{\z0U_k}$, so that in particular we have the estimate
below,  which relies upon \e.u8/.	
 	\md0
\norm b^{\z0U_k}_{\text{trun}}.p.\seq{}\ze2^{k/100},\qquad
1<p<\zI.
\emd
(By a suitable definition of the strong maximal function $M$, one
can deduce this inequality   from the $L^p$ bounds for
$M$.)
 We apply this inequality far away from the set $U$. For the set
$W_\zl=\ZR^2-\bigcup_{R\in \z0U_k}\zl R$, $\zl>1$, we have the
inequality 	
\md1\Label  e.trun
 \norm
b^{\z0U_k}_{\text{trun}}.L^p(W).\seq{}\ze2^{k/100}\zl^{-100},
\qquad 1<p<\zI.
 \emd

\bigskip 

We shall need a refined decomposition of the collection $\z0W$, the motivation
for which is the following calculation.  Let $\z0W'\subset\z0W$.  For 
$n=(n_1,n_2)\in\ZZ^2$, set $\z0W'(n):=\{R'\in \z0W'\mid 
\abs{R'_j}=2^{n_j},\ j=1,2\}$.  In addition, let 
\md0
B(\z0W',n):=\sum_{R'\in\z0W'(n)}\sum_{\substack{R\in\z0U_k
\\ R'< R}}{\ip b,v_{R'}.\z8{\ip b,v_R.}
v_{R'}\z8{v_R}}.
\emd
And set $B(\z0W')=\sum_{n\in\ZZ^2}B(\z0W',n)$. 

Then, in view of \e.doubleor/, we see that $B(\z0W',n)$ and $B(\z0W',n')$ are
orthogonal if $n$ and $n'$ differ by at least $3$ in either coordinate.  
Thus, 
\md0
\Bigl\lVert \sum_{n\in\ZZ^2}B(\z0W',n)\Bigr\rVert_2^2
{}\leq{}3
\sum_{n\in\ZZ^2} \lVert  B(\z0W',n) \rVert_2^2.
\emd

The rectangles $R'\in\z0W(n)$ are all translates of one another.  Thus, taking
advantage of the rapid spatial decay of the  wavelets, we can estimate 
\md0
 \lVert  B(\z0W',n) \rVert_2^2
{}\seq{}
\sum_{R'\in\z0W(n)}\int \ABs{ \frac{\abs{\ip
b,v_{R'}.}}{\sqrt{\abs {R'}}}(\zq_{R'}*1_{R'})
b^{\z0U_k}_{\text{trun}}}^2 \; dx \emd In this display, we let
$\zq(x_1,x_2)=(1+x_1^2+x_2^2)^{-10}$ and for rectangles $R$,
$\zq_R(x_1,x_2)=\zq(x_1\abs{R_1}^{-1},x_2\abs{R_2}^{-1})$.  Note
that $\zq_R$ depends only on the dimensions of $R$ and not its
location.

Continuing, note the trivial inequality 
 $\int (\zq_R*f)^2g\;dx\seq{}\int\abs{f}^2\zq_R*g\;dx$.
We can estimate   
\md3 \nonumber 
     \lVert B(\z0W')\rVert_2^2
     {}\seq{}& 
     \sum_{  R'\in\z0W'} \abs{ \ip b,v_{R'}.}^2
\Bigl\{\abs{R'}^{-1}
\int_{R'}M(\abs{b^{\z0U_k}_{\text{trun}}}^2)\; dx\Bigr\}  
\\{}\seq{}& 
\ABs{\bigcup_{R'\in\z0W'}R'} \sup_{R'\in\z0W'} \text{avg}(R').
 \Label e.repeat
\emd
Here we take $\text{avg}(R'):=\abs{R'}^{-1}
\int_{R'}M(\abs{b^{\z0U_k}_{\text{trun}}}^2)$. 

The terms $\text{avg}(R')$ are essentially of the order of magnitude 
$\ze^2$ times a the scaled distance between $R'$ and the open set $U$. 
To make this precise requires a decomposition of the collection $\z0W$. 

\bigskip

For integers $l>k$ and $m\ge0$, set $\z0W(l,m)$ to be those $R'\in\z0W$ which
satisfy these three conditions.  
\begin{itemize}
\item 
First, $\text{avg}(R')\le\ze2^{-4l}$ if $m=0$ and $\ze2^{-4l+m-1}
<\text{avg}(R')\le\ze2^{-4l+m} $ if $m>0$.
\item 
 Second, there is
an $R\in\z0U_k$ with $R'<R$  and  $R'\subset 2^{l+1}R$.
\item
  Third, for
every $R\in\z0U_k$ with $R'<R$, we have $R'\not\subset
2^{l+1}R$.
Certainly, this collection of rectangles is empty if $l\le k$.   
\end{itemize}

We   see that 
\md0
\ABs{\bigcup_{R'\in\z0W(l,m)}R'}\seq{}\min(2^{2lp}, 2^{-mp/2}),
\qquad 1<p<\zI. \emd
The first estimate follows since the rectangles $R'\in\z0W(l,m)$ are contained
in the set $\{M1_U\ge2^{-2l-1}\}$.  The second estimate follows from \e.trun/.

But then from \e.repeat/ we see that for $m>0$,
\md0
\lVert B(\z0W(l,m))\rVert_2^2\seq{}\ze2^{-4l+m}\min(2^{2lp},2^{-mp/2})
{}\seq{}\ze2^{-(m+l)/10}.
\emd
In the case that $m=0$, we have the bound $2^{2lp}$.
This  is obtained by taking the minimum to be $2^{2lp}$ for
$p=5/4$ and $0<m<\frac{11}8l$.  For $m\ge\frac{11}8l$ take the
minimum to be $2^{-mp/2}$ with $p=4$.

This last estimate is summable over $0<k<l$ and $0<m$ to at most ${}\seq{}\ze$, and so
completes this case.

 %%%%%%%%%%%%%%%%%%%%%%%%%    <_1
\subsubsection{The Partial Orders `$<_j$', $j=1,2$.}

 We treat the case of $R'<_1R$, while  the case of
  $R'<_2R$ is same by symmetry.  The structure of this partial order
  provides some orthogonality
  in the first variable, leaving none in the second variable.
  Bounds for  the expressions from
  the second variable are   derived from    a cognate of a
Carleson measure estimate.

  There is a basic calculation that we perform for a subset $\z0W'\subset \z0W$.  For
  an integer $n'\in\ZZ$ define $\z0W'(n'):=\{R'\in\z0W'\mid
\abs{R'_1}=2^{n'}\}$, and   \md0
B(\z0W',n'):=\sum_{R'\in
\z0W'(n')}\sum_{\substack{R\in\z0U_k\\R'<_1R}
}   \ip b,v_{R'}.\z8{\ip b,v_R.} v_{R'}\z8{v_R}.
\emd
Recalling \e.doubleor/, if  $n'$ and $n''$ differ by more than
$3$,   then $B(\z0W',n')$ and $B(\z0W',n'')$ are orthogonal.

  Observe that for $R'$ and $R$ as in the sum defining $B(\z0W',n)$, we have the
  estimate
  \md1  \Label e.v<z
  \abs{v_{R'}(x)\z8{v_R}(x)}\seq{}
  (\abs R\abs{R'})^{-1/2}
  \text{dist}(R',R)^{1000}\zq_{R'}*1_{R'}(x),\qquad
  x\in\ZR^2.
  \emd
  In this display, we are using the same notation as before, 
   $\zq(x_1,x_2)=(1+x_1^2+x_2^2)^{-10}$ and for rectangles
$R$, $
\zq_R(x_1,x_2)=\zq(x_1\abs{R_1}^{-1},x_2\abs{R_2}^{-1})$. 
  In addition, $\text{dist}(R',R):=M1_R(c(R'))$, with $c(R')$ being the
  center of $R'$. [This ``distance" is more properly the inverse
of a distance that takes into account the scale of the rectangle
$R$.]

  \smallskip

  Now define
  \md1 \Label e.zbdef
  \zb(R'):=\sum_{\substack{R\in\z0U\\
R'<_1R}}\abs{R}^{-1/2}\abs{\ip b,v_R.} \text{dist}(R',R)^{1000}   .
\emd
  The main point of these observations and definitions is this.  For the
  function $B(\z0W'):=\sum_{n'\in\ZZ}B(\z0W',n')$, we have
  \md2
  \lVert B\rVert_2^2\seq{}&
  \sum_{n'\in\ZZ}\lVert B(\z0W',n')\rVert_2^2
  \\{}\seq{}&
   \sum_{n'\in\ZZ}\int_{\ZR^2}
  \Biggl[ \sum_{R'\in\z0W'(n')}  \abs{\ip
b,v_{R'}.}\zb(R')\abs{R'}^{-1/2}   \zq_{R'}*1_{R'}\Biggr]^2\; dx
  \\{}\seq{}&    \sum_{n'\in\ZZ}
   \int_{\ZR^2}
  \Biggl[ \sum_{R'\in\z0W'(n')}  \abs{\ip
b,v_{R'}.}\zb(R')\abs{R'}^{-1/2}    1_{R'}\Biggr]^2\; dx
    \emd
 At this point, it occurs to one to appeal to the Carleson
measure property associated to the coefficients $\abs{\ip
b,v_{R'}.}\abs{R'}^{-1/2}$. This necessitates that one proves
that the coefficients  $\zb(R')$ satisfy a similar
condition, which doesn't seem to be true in general.  A slightly weaker
condition is however true.

\bigskip

To get around this difficulty, we make a further diagonalization
of the terms $\zb(R')$ above.
 For integers $\zn\ge\zn_0$, $\zm\ge1$ and a
rectangle $R'\in\z0W$, consider rectangles $R\in\z0U_k$ such
that \md0
R'<_1R,\qquad 2^{-v}\le\dist(R',R)\le{}2^{-v+1}
,\qquad 2^\zm\abs{R'}=\abs{R}.
\emd
[The quantity $v_0$ depends upon the particular subcollection
$\z0W'$ we are considering.]
We denote one of these rectangles as  $\zp(R')$.

An important geometrical fact is this.  We have $\zp(R')\subset
2^{v+\zm+10}R'_1\times 2^{v+10}R_2'$.  And in particular, this
last rectangle has measure ${}\seq2^{2v+\zm}\abs{R'}$.

Therefore, there are at most $O(2^{2v})$ possible choices for $\zp(R')$.
[Small integral powers of $2^v$ are completely harmless because of the 
large power of $\text{dist}(R',R)$  that appears in \e.zbdef/.]

Our purpose is to bound this next expression by a term which
includes a power of $\ze$, a small power of $2^v$ and a power of
$2^{-\zm}$.  Define
\md2
S(\z0W',\zn,\zm):={}& \sum_{n'\in\ZZ}\int_{\ZR^2}
  \Biggl[ \sum_{R'\in\z0W'(n')}  \frac{\abs{\ip
b,v_{R'}. \ip b,v_{\zp(R')}.
} }{\sqrt{\abs{R'}\abs{\zp(R')}}}
\zq_{R'}*1_{R'}\Biggl]^2\; dx
\\{}\seq{}&
  \sum_{n'\in\ZZ}\int_{\ZR^2}
  \Biggl[ \sum_{R'\in\z0W'(n')}  \frac{\abs{\ip
b,v_{R'}.\ip b,v_{\zp(R')}.} }{\sqrt{\abs{R'}\abs{\zp(R')}}}
1_{R'}\Biggl]^2\; dx
\\{}={}&     \sum_{n'\in\ZZ}  \sum_{R'\in\z0W'(n')}
 \frac{\abs{\ip
b,v_{R'}.\ip b,v_{\zp(R')}.}} {\sqrt{\abs{R'}\abs{\zp(R')}} }
\sum_{\substack{R''\in\z0W'(n')\\ R''\subset R'} }
\sqrt{\frac{\abs{R''}}{\abs{\zp(R'')}}}  \abs{\ip
b,v_{R''}.\ip b,v_{\zp(R'')}.}
\emd

The innermost  sum can be bounded this way.  First $\norm
b.BMO(rec).\le\ze$, so that
\md0
  \sum_{ R''\subset R'} \abs{\ip
b,v_{R''}.}^2  {}\le{}\ze^2\abs{R'}.
\emd
Second, by our geometrical observation about $\zp(R')$,
  \md0
  \sum_{ R''\subset R'}\frac{\abs{R''}}{\abs{\zp(R'')}} \abs{\ip
b,v_{\zp(R'')}.}^2
{}\seq{}\ze^22^{2v}\abs{R'}.
\emd
In particular, the factor $2^u$ does not enter into this estimate. 

This means that
\md2
S(\z0W,\zn,\zm)\seq{} & \ze^22^{2v}
\sum_{R'\in\z0W'} {\sqrt{\frac{\abs{R'}}{\abs{\zp(R')}}}}
\abs{\ip b,v_{R'}.\ip b,v_{\zp(R')}.}
\\{}\seq{}&
  \ze^22^{2v-\zm/2}
  \Bigl[ \sum_{R'\in\z0W'}\abs{ \ip b,v_{R'}.}^2
  \sum_{R\in\z0U_k}\abs{\ip b,v_R.}^2  \Bigr]^{1/2}
 \\{}\seq{}&
  \ze^22^{2v-\zm/2}\ABs{\bigcup_{R'\in\z0W' }R'}^{1/2}
  \emd

\bigskip

The point of these computations is that a further trivial
application of the Cauchy--Schwartz inequality proves that
 \md0
\norm B(\z0W').2.\seq{} \ze 2^{-100\zn_0}
\ABs{\bigcup_{R'\in\z0W' }R'}^{1/4}
 \emd
where $\zn_0$ is the largest integer such that for all
$R'\in\z0W'$ and $R\in\z0U_k$, we have
$\dist(R',R)\le2^{-\zn_0}$.

   \smallskip
   We shall   complete this section by decomposing $\z0W$ into
subcollections for which this last estimate
summable to  $\ze2^{-k}$.
Indeed, take $\z0W_v$ to be those $R'\in\z0W$ with $R'\not\subset
2^vR$ for all    $R\in\z0U_k$ with $R'<_1R$.  And there is an
$R\in\z0U_k$ with    $R'\subset2^{v+1} R$  and $R'<_1R$.  Certainly, we 
need only consider $v\ge{}k$.  

It is clear that this decomposition of $\z0W$ will conclude the
treatment of this partial order.

\subsubsection{The partial order `$\simeq$'}
 
We now consider the case of $R'\simeq R$, which is less 
subtle as there is no orthogonality to exploit
and the Carleson measure estimates are more directly applicible.
We prove the bound
\md0 \NORm\sum_{R'\in\z0W}\sum_{\substack{R\in\z0U
\\R'\simeq R}}{\ip b,v_{R'}.\z8{\ip b,v_R.} 
v_{R'}\z8{v_R}}.2.\seq{} K_{\zd}\epsilon^{1/3} .
\emd
 
The diagonalization in space takes two different forms.  
  For $\zl\ge2^k$ and
$R\in\z0U_k$ set $\zs(\zl,R)$ to be a choice of $R'\in\z0W$ with 
$R'\simeq R$ and $R'\subset 2\zl R$. (The 
definition is vacuous for $\zl<2^k$.) 
 It is clear 
that we need only consider ${}\simeq\zl^2$ choices of these 
functions $\zs(\zl,\cdot)\mid\z0U_k\longrightarrow\z0W$.  There is 
an $L^1$ estimate which allows one to take advantage of 
the spatial separation between $R$ and $\zs(\zl,R)$. 
  \md2 
\NORm\sum_{{R\in\z0U_k}}{\ip b,v_{\zs(\zl,R)}.\z8{\ip b,v_R.} 
v_{\zs(\zl,R)}\z8{v_R}}.1.\seq{} & 
        \zl^{-100}\sum_{R\in\z0U_k}\abs{ 
        \ip b,v_{\zs(\zl,R)}.\z8{\ip b,v_R.}} 
\\{}\seq{}& 
     \zl^{-100}\Bigl[   \sum_{R\in\z0U_k}\abs{ 
         \ip b,v_{\zs(\zl,R)}.}^2  \sum_{R\in\z0U_k}\abs{ 
         \ip b,v_{R}.}^2\Bigr]^{1/2} 
\\{}\leq{}& K_{\zd}{\epsilon} 
     \zl^{-90}. 
\emd This estimate uses \e.u8/ and is a very small estimate. 
 
To complete this case we need to provide an estimate in $L^4$. 
Here, we can be quite inefficient.  By 
Cauchy--Schwartz and the  Littlewood--Paley inequalities, 
 \md2 
\NORm\sum_{{R\in\z0U}} {\ip b,v_{\zs(\zl,R)}.\z8{\ip b,v_R.} 
v_{\zs(\zl,R)}\z8{v_R}}.4.\seq{} & 
    \NORm\Bigl[\sum_{{R\in\z0U}}\abs{{\ip b,v_{\zs(\zl,R)}. 
v_{\zs(\zl,R)}}}^2\Big]^{1/2} .4. 
\NORm\Bigl[\sum_{{R\in\z0U}}\abs{{\ip 
b,v_R.}\z8{v_R}}^2\Bigr]^{1/2}.4. \\{}\seq{}& \zl. \emd This follows 
directly from the $\bmo$ assumption on $b$. Our proof is complete.

%%%%%%%%%%%%%%%%%%%%%%%%%%%%%%%%%%%% 
 
\appendix 
 
\section{A Remark on Journ\'e's Lemma} 
 
Let $U$ be an open set of finite measure in the plane.  Let 
$\z0R^*(U)$ be those maximal dyadic rectangles in $\z0R$ that are 
contained in $U$.   Define   for each $0<\zd<1$ and $i,j\in\{1,2\}$, 
  \md0
 V_{\zd,i, j}=\{M_i\ind {\{M_j 
\ind{}U>\zd\}}>\zd\} , 
  \emd
  and $V_\zd=V_{\zd,1, 2}\cup V_{\zd,2, 1}$.   For each $R\in\z0R(U)$ 
set \md0 \zm_\zd(R)=\sup\{\zm>0\mid \zm R\subset V_\zd\}. \emd The 
form of Journ\'e's Lemma we need is 

\bthm l.journe For each 
$0<\zd,\ze<1$ and each open set $U$ in the plane of finite 
measure, \md0 \sum_{R\in\z0R^*(U)} 
\zm_\zd(R)^{-\ze}\abs{R}\seq\abs U. \emd \ethm 
 
Journ\'e's Lemma is the central tool in verifying the Carleson 
measure condition, and points to the central problem in two 
dimensions:  That there can be many rectangles close to the 
boundary of an open set. 
 
Among the references we could find in the literature 
\cite{journe,j1}, the form of Journ\'e's Lemma cited and proved is 
relative to a strictly larger quantity than the one we use, 
$\zm_\zd(R)$ above. To define it, for any rectangle $R$, denote it 
as a product of intervals $R_1\times R_2$. 
 Set $M$ to be the strong maximal function.
 Then, for open set $U$ 
of finite measure and $R\in\z0R(U)$, set (taking $\zd=1/2$ for 
simplicity) \md0 \zn(R)=\sup\{\zn>0\mid \zn R_1\times R_2\subset 
\{M\ind U >{1/2}\}\}. 
 \emd 
Thus, one only measures how deeply $R$ is in the enlarged set 
in one direction. The Lemma above then holds for $\zn(R)$, with 
however a slightly sharper form of the sum than we prove here. 
 
In addition, note that one measures the depth of $R$ with respect 
to  a simpler set, $\{M\ind U>1/2\}$.  We did not use this 
simplification in our proof as   the strong 
maximal function $M$ does not act boundedly on $L^1$ of the plane. 
 
 There are however examples which show that that the quantity 
$\zn(R)$ can be much larger than $\zm(R)$.   Indeed, consider a 
horizontal row of evenly spaced squares.  For a square $R$ in the 
middle of this row,  $\zn(R)$ will be quite big, while $\zm(R)$ 
will be about $1$ for all $R$.  Thus we give a proof of our form 
of Journ\'e's Lemma. 
 
 \begin{proof}[Proof of \l.journe/.]   We can assume that $1/2\le\zd<1$ 
as the 
terms $\zm_\zd(R)$  decrease as $\zd$ increases. Fix $\zm\ge1$. 
Set $\z0S$ to be those rectangles in $\z0R^*(U)$ with $ 
\zm\le\zm_\zd(R)\le2\zm$. It suffices to show that \md0 
 \sum_{R\in\z0S}\abs{R}\seq(1+\log \zm)^2\abs U. 
\emd 
 For then this estimate is summed over $\zm\in\{2^k\mid k\in\ZZ\}$. 
 
 In showing this estimate, we can further assume that for all 
$R,R'\in\z0S$, 
writing $R=R_1\times R_2$ and likewise for $R'$, that if for 
$j=1,2$, $\abs{R_j}>\abs{R'_j}$ then 
$\abs{R_j}>16\zm(1-\zg)^{-1}\abs{R'_j}$, where we set 
$\zg=\zd^{1/3}$. This is done by restricting $\log_2\abs{R_j}$ to 
be in an arithmetic progression of difference ${}\simeq\log \zm$.  This 
neccessitates the division of all rectangles into ${}\seq(1+\log \zm)^2$ 
subclasses  and so we 
prove the bound above without the logarithmic term.

We define a ``bad" class of rectangles  $\z0B=\z0B(\z0S)$
as follows.  For $j=1,2$, let $\z0B_j(\z0S)$ be 
 those rectangles $R$ for which there are rectangles 
 \md0 
 R^1,R^2,\ldots,R^K\in\z0S-\{R\}, 
 \emd 
 so that for each $1\le{}k\le{}K$, $\abs{R^k_j}>\abs{R_j}$, and 
 \md0 
 \ABS{ R\cap \bigcup_{k=1}^K R^k}>\zg\abs R. 
 \emd 
 Thus $R\in\z0B_j$ if it is nearly completely covered in  the $j$th direction of 
the plane.  
Set $\z0B(\z0S)=\z0B_1(\z0S)\cup\z0B_2(\z0S)$.      It follows that 
if $R\not\in\z0B(\z0S)$, it is not covered in both the vertical and horizontal 
directions, hence 
\md0 \ABS{ R\cap 
\bigcap_{R'\in\z0S-\{R\}}(R')^c}\ge(1-\zg)^2\abs R. 
\emd 
And, since all $R\subset U$, it follows that \md0 
\sum_{R\in\z0S-\z0B}\abs R\le(1-\zg)^{-2}\abs U. \emd 

\medskip 
 
Thus, it remains to consider seperately  the set of rectangles $\z0B_1(\z0S)$ 
and $\z0B_2(\z0S)$.      Observe that for any collection $\z0S'$, 
$\z0B_j(\z0S')\subset\z0S'$ as follows 
immediately from the definition.  Hence 
$\z0B_1(\z0B_2(\z0B_1(\z0S)))\subset{} \z0B_1(\z0B_1(\z0S))$.  And we argue that this last set is empty.  
As our definition of $V_\zd$ and $\zm(R)$ is symmetric with respect to the 
coordinate axes, this is enough to finish the proof. 
 
 We argue that $      \z0B_1(\z0B_1(\z0S))$ is empty by contradiction. 
Assume that $R\in\z0B'$.               Consider 
  those rectangles $R'$ in $\z0B_1(\z0S)
  $ for which $(i)$ $\abs{R'_1}>\abs{R_1}$  and 
$(ii)$ 
$R'\cap R\not=\emptyset$. 
Then \md0 \ABS{R\cap\bigcup_{R'\in\z0C}R'}\ge\zg\abs R. 
\emd 
Fix a  one of these rectangles $R'$ with 
 $\abs{R'_1}$ being minimal.  We then claim that 
$8\zm R'\subset \{\z0M_1\ind{\{\z0M_2\ind{} U>\zd\}}>\zd\}$, which 
contradicts the assumption that $\zm(R')$ is no more than $2\zm$. 
 
Indeed,  all the rectangles in $\z0B_1(\z0S)$ are themselves covered in 
the first  coordinate axis. We see that the the set  $\{\z0M_2\ind U>\zg^3\} $ 
contains the rectangle $R''_1\times\zg^{-1} R_2$, in which $R_2$ 
is the second coordinate interval for the rectangle $R$ and 
$R''_1$ is the dyadic interval that contains $R'_1$ and has 
measure $8\zm(1-\zg)^{-1}\abs{R_1'}\le{} 
\abs{R''_1}<16\zm(1-\zg)^{-1}\abs{R_1'}$. 
 
But then the rectangle $\zg^{-1}R''_1\times \zg^{-1}R_2$  is 
contained in  $\{\z0M_1\ind{ \{\z0M_2\ind U>\zg^3\}}>\zg^3\}$. And 
since $8\zm R'$ is contained in this last rectangle, we have 
contradicted the assumption that $\zm(R')<2\zm$. 
 
\end{proof} 
\endgroup   %  forthe spacing of the paragraphs 
 
%%%%%%%%%%%%%%%%%%%%%%%%%%%%%%%%%%%% 

\bigskip
{\parindent=0pt\baselineskip=12pt\obeylines
Sarah H. Ferguson
 Department of Mathematics
  Wayne State University
Detroit MI 48202
\smallskip
\tt sarah@math.wayne.edu \tt http://www.math.wayne.edu/\~{}sarah }
 \medskip
 {\parindent=0pt\baselineskip=12pt\obeylines 
Michael T. Lacey
School of Mathematics
Georgia Institute of
Technology Atlanta GA 30332 
\smallskip 
\tt lacey@math.gatech.edu \tt http://www.math.gatech.edu/\~{}lacey } 
\medskip


\begin{thebibliography}{99} 
 
\bibitem{cf} S.--Y. A. Chang  and R. A. Fefferman.  {``Some recent 
developments in 
Fourier analysis and $H\sp p$-theory on product domains"} {\it
Bull. Amer. Math. Soc. (N.S.) } {\bf 12} (1985)  1--43 



\bibitem{one} R.R.~Coifman,  R.~Rochberg, G.~Weiss. {``Factorization 
theorems for Hardy spaces in several variables."} Ann. of Math. (2) {\bf
103} (1976),
 611--635.
 
\bibitem{f} R. A. Fefferman. {``Bounded mean oscillation on the polydisk"} 
Ann. Math. {\bf 110} (1979) 395-406 


\bibitem{fs} S. H. Ferguson and C. Sadosky. {``Characterizations of 
bounded mean oscillation on the polydisk in terms of Hankel 
operators and Carleson measures''} {\it J. 
D'Analyse Math.} {\bf 81} (2000) 239-267 
 
\bibitem{gs} R. Gundy and E. Stein. {``H$^p$ theory for the
poly-disc"} {\it Proc. Nat. Acad. Sci.} {\bf 76} (1979) 1026-1029
 
\bibitem{journe} J.--L. Journ\'e {``A covering lemma for product 
space."} {\it 
Proc. Amer. Math. Soc.} {\bf 96} 593---598. MR87g:42028 
 
\bibitem{m}  Yves Meyer. f{\sl 
 Wavelets and operators.}
 Translated from the 1990 French original by D. H. 
Salinger. Cambridge Studies in Advanced Mathematics, 37. Cambridge University 
Press, Cambridge, 1992. xvi+224 pp. ISBN: 0-521-42000-8 MR94f:42001.

\bibitem{n} Z. Nehari.   Ann. of Math. (2) {\bf 65} (1957),
 153-162. MR21\#7399.
 
\bibitem{j1} J. Pipher. ``Journ\'e's covering lemma and its extension to 
higher 
dimensions" {\it Duke Math. J.} {\bf 53} (1986), no.~3, 683--690;
MR 88a:42019 
 
\bibitem{sw} E. M. Stein and G. Weiss.{\sl 
Introduction to Fourier Analysis on 
Euclidean Spaces.}
 Princeton Univ. Press Princeton, 1971 
\end{thebibliography}
\end{document}